\documentclass[10pt]{amsart}
\usepackage{indentfirst}
\usepackage{latexsym}
\usepackage{amssymb}
\usepackage{amsmath}
\usepackage{mathrsfs}
\usepackage{graphicx}

\begin{document}

\title{ALEXANDROV PAR EXCELLENCE}
\author{S.~S. Kutateladze}
\address[]{
Sobolev Institute of Mathematics\newline
\indent 4 Koptyug Avenue\newline
\indent Novosibirsk, 630090\newline
\indent Russia}
\email{
sskut@member.ams.org
}
\begin{abstract}
This is a short tribute to Alexandr Alexandrov
(1912--1999).
\end{abstract}
\date{January 7, 2007}

\maketitle

{\baselineskip=.95\baselineskip

Many mass media declared Grisha Perelman
the 2006 Scientist of the Year  in recognition of the Fields Medal award for
his proof of the Poincar\'e conjecture.
Perelman is  the last postgraduate student supervised by
Alexandr Danilovich Alexandrov (1912--1999) the anniversary of whose birth
we celebrate this year.

Many years ago some visiting professor persuaded Alexandrov in my presence
to make a demarche  at the public maintenance of a mediocre  geometric thesis.
The defendant was a Jew, and the final argument of the Muscovite professor
was xenophobia. He pointed to the picturesque and rather gloomy figure
of a complete stranger to the event at the end of a corridor and began blaming
us in the sense that if  we  disagreed to his proposal
then ``everyone would be like that down here.''
By the vicissitudes of fate, the lad in the corner  was  Perelman,
a postgraduate of Alexandrov, a student of Professor  Yuri\u\i{} Bugaro who 
was himself
a student of Alexandrov.

The love and hatred to Alexandrov stemmed from the same sources.
His reviews and opinions were welcome and appreciated, but his approaches and areas
of research were silenced if not scorned. He was accused of zionism, but many bet
and counted  upon his antisemitism. His communistic beliefs were  blasphemed obscenely,
but he was humbly requested to write a letter or two to the Central Committee
of the Communist Party of the USSR or  {\it The Communist} magazine.
His philosophical essays were spit upon furiously, but the same despisers
required that their students used Alexandrov's popular writings at the final
examinations in philosophy which were obligatory for admittance to the public maintenance of  theses.
The professorship of St. Petersburg is full of raptures about the
palace, fountain, and park ensemble of Peterhof, but most of Alexandrov's colleagues
will never forgive the sage decision of Rector Alexandrov
who suggested to build a new university campus in Peterhof.
During the years of Gorbi's ``perestroika'' Alexandrov was accused in confessing
lysenkoism but decorated with the Order of the  Red Banner of Labor
for his efforts in safegarding and propelling genetics and selection
in the USSR.
So were the scales of  Alexandrov's personality.

Alexandrov often remarked that a man is what he does.
Alexandrov did the business of geometry.
It is more appropriate to speak about geometry as
his most adorable part of the universal science, mathematics.
Mac Lane, a cofather of category theory,
had coined the term ``working mathematician''
which is rather close to the routine
``workman'' or even ``toiler.''
Mac Lane wrote not only about  work in mathematics but also about
criteria of excellence in mathematics.  Excellent mathematics by Mac Lane
should be inevitable, illuminating, deep, relevant, responsive, and timely.
Excellent mathematics is exercised by excellent mathematicians,
mathematicians  {\it par excellence}. Such was Alexandrov.

Alexandrov contributed to mathematics under the slogan:
``Retreat to Euclid.'' He remarked that
``the pathos of contemporary mathematics is the return to Ancient Greece.''
Hermann Minkowski revolutionized the theory of numbers
with  the aid of the synthetic geometry of convex figures.
The ideas and techniques of the geometry of numbers comprised
the fundamentals od functional analysis which was created by Banach.
The pioneering studies of Alexandrov continued the efforts
of Minkowski and enriched geometry with the methods
of measure theory and functional analysis.
Alexandrov accomplished the turnround to the ancient synthetic geometry
in a much  deeper and subtler sense than it is generally acknowledged today.
Geometry in the large  reduces in no way to overcoming
the local restrictions of differential geometry which bases
upon the infinitesimal methods and ideas of Newton,
Leibniz, and Gauss.

The works of Alexandrov made tremendous progress in the theory of
mixed volumes of convex figures. He proved  some fundamental
theorems on convex polyhedra that are celebrated alongside
the theorems of Euler and Cauchy.
While discovering a~solution of the Weyl problem,
Alexandrov suggested  a new synthetic method for proving the theorems of existence.
The results of this research ranked the name of Alexandrov  alongside
the names of Euclid and Cauchy.

Another outstanding contribution of Alexandrov to science is the creation
of the intrinsic geometry of irregular surfaces. He suggested his
amazingly visual and powerful method of cutting and gluing.
This method enabled him to solve many extremal problems of the
theory of manifolds of bounded curvature.

Alexandrov developed the theory of metric spaces with one-sided constraints on curvature.
This gave rise to the class of metric spaces generalizing the Riemann spaces in the sense
that these spaces are furnished with some curvature, the basic concept
of Riemannian geometry.
The research of Alexandrov into the theory of manifolds with bounded curvature
prolongates and continues the traditions of Gauss, Lobachevski\u\i{}, Poincar\'e, and Cartan.

Alexandrov enriched the methods of differential geometry by
the tools of functional analysis and measure theory,
driving mathematics to its universal status of the epoch of Euclid.
The mathematics of the ancients was geometry
(there were no other instances of mathematics at all).
Synthesizing geometry with the remaining areas of the today's mathematics,
Alexandrov climbed to the antique ideal of the universal science incarnated in
mathematics.
Return to the synthetic methods of {\it mathesis universalis\/}
was inevitable and unavoidable  as illustrated in geometry with the beautiful
results of the disciples and followers of Alexandrov like
Gromov, Perelman, Pogorelov, and Reshetnyak.

The first and foremost Russian geometer of the nineteenth century was
Nikola\u\i{} Ivanovich Lobachevski\u\i.

Alexander Danilovich Alexandrov became
the first and foremost Russian geometer of the twentieth century.

}

\end{document}